\documentclass[runningheads,a4paper]{llncs}

\usepackage{amssymb}
\usepackage{graphicx}

\usepackage{url}
\urldef{\mailsa}\path| mahyunst@yahoo.com|
\newcommand{\keywords}[1]{\par\addvspace\baselineskip
\noindent\keywordname\enspace\ignorespaces#1}

\begin{document}

\mainmatter  

\title{The Ontology of Knowledge Based Optimization}

\titlerunning{The Ontology of Knowledge Based Optimization}
 
\author{Mahyuddin K. M. Nasution}

\authorrunning{M. K. M. Nasution}
\institute{Mathematic Department, Fakultas Matematika dan Ilmu Pengetahuan Alam\\
Universitas Sumatera Utara, Padang Bulan 20150 USU Medan Indonesia\\
\mailsa\\}

\toctitle{SIMANTAP 2010}
\tocauthor{}
\maketitle

\begin{abstract}
Optimization has been becoming a central of studies in mathematic and has many areas with different applications. However, many themes of optimization came from different area have not ties closing to origin concepts. This paper is to address some variants of optimization problems using ontology in order to building basic of knowledge about optimization, and then using it to enhance strategy to achieve knowledge based optimization.
\end{abstract}

\section{Introduction}
The optimization refers to choosing the best element from some set of available alternatives \cite{diwekar2008}. Simply, the optimization related to minimizing or maximizing a function by systematically choosing the values in/from/to real. In some respects the optimization embedded to economic, efficiency, and effectiveness. Therefore, in simple, a optimization problem can be defined as a function 
\begin{equation}
\label{pers:utama}
f : X\rightarrow {\bf R}
\end{equation}
where $\exists x_0\in X, \forall x_i\in X$ such that a minimizaton is $f(x_0)\leq f(x_i)$ or a maximization $f(x_0)\geq f(x_i)$ \cite{guler2010}. In general, a set of variables $X$ is a collection of variables representing attributes of entities as instances in studies which need an optimization. Some of them are the standard variables (static), but do not few dynamically as options. However, in many researches and studies, or the research behaviors about optimization, those variables on fixed condition is either in setting or in concept, such as many examples in course books of optimization, such that the optimization formulation proved to be stiff, and indeed do not applicable.

Most optimization-related research in the airline industry, for example, are about areas: network design and schedule construction; fleet assignment; aircraft routing; crew scheduling; revenue management; irregular operations; air traffic control and ground delay program \cite{yu2008}. Which of all consider the green concept \cite{buttazzo2009}, for example. 
Therefore, each time the optimization problem formulated, but it is not optimization. The reasons, many studies consider the properties of variables setting, but they are not sure the characteristics of variables, so the latent variable can not be disclosed, need to both external and internal considerations as well as for providing impact in the optimization variables, such as the reliability and trusty. Therefore, this paper has a task to describe some problems of optimization using ontology, and to find out ties between problems based on adaptation concept. This paper divided into four sections, one remainder section is about the concept and some problem statements, and then they are expressed in Discussion section.

\section{The Concepts and Some Open Problem Statements}

The research in optimization area are becoming more flamboyant with the emergence of new area of studies involving resources (implicitly \cite{fisch2003}). Such studies include the exploration of new set theory such as fuzzy set \cite{jensen2005}. The optimization involves many themes and interest of every branch of mathematics as knowledge that exist today, which has the task in accordance with its purpose to identify the new things around optimization and assessment of object relations with other agencies, which also involves the main agenda of life such as environment, human rights, social welfare, and justice, where the green computing and nano-technology have played their roles\cite{audet2005,bock2005,baker2007,alves2008,chinchuluun2008,kugler2008,kosmidou2008,lowen2008,chaovalitwongse2009,papajorgji2009,gonzalez2010}. These do not only experienced a shift in accordance with the challenges and problems faced to get optimization, but also caused by changes how to approach the issues by the scientists. 

This shift is a response to the pressures of internal and external. Internal pressure is a jigsaw puzzle that is still hidden and not answered within the given optimization paradigm, caused by the approaches or methodologies that have not been able to explain  phenomena or any event which can be observed, while the external pressure comes from the style of thought and the flow of thought. The flow of thought such as ontology is used to find a red thread between the paradigm and phenomena \cite{pickard2007}, thus optimization studies can be clearly disclosed in both research opportunities and applications. Existing red threads include heuristic rules, inductive, and deductive. However, it is also not effective to directly employ to many phenomena, this is because 
\begin{enumerate}
\item the heuristic rule requires the scientist to define a specific rule for each specific type of optimization problem, thus that is not adaptive for different situations {\bf [Adapdation 1]},
\item the inductivity trains a scientist model individually and cannot be adapted to the another {\bf [Adapdation 2]}, and
\item deductivity can deal with different instances simultaneously, but it cannot make use of the prediction {\bf [Adapdation 3]}.
\end{enumerate}

One of optimization phenomena is about feasible solution (is called optimal solution) that minimizes or maximizes the objective function, where the feasible region and the objective function present convexity or not. If (\ref{pers:utama}) as objective function (OF), cost function or energy function does not present convexity, then there may be several {\it local optima} (minima or maxima), where local optima $x^*$ is defined as a point for which there exists some $\delta>0$ so that for all $x$ such that $\|x-x^*\| ~(\leq\vee \geq)~ \delta$; the expression $f(x^*) ~(\leq\vee\geq)~ f(x)$ holds, whether the opposite is a global optima or not. Based on this conditions, there are two categories of optimization which always use the programming term as an emphasis on the use of iteration in complexity to get the solution of problem. First category is convex programming, i.e., optimization studies when the objective function (\ref{pers:utama}) is convex and the constraints, if any, also form a convex set. In general, the convex programming as follow \cite{burachik2008}.

\begin{definition} [Convex optima]
\label{def:ko}
 Let a real vector space $X$ together with a convex, real-valued function 
\begin{equation}
\label{pers:convexoptima}
f:{\cal X}\rightarrow{\bf R}
\end{equation}
such that $f({\bf x}^*) \leq f({\bf x})$, $\forall {\bf x}\in X$, $\exists {\bf x}^*\in{\cal X}$ for the number $f({\bf x})$ is smallest, 
where ${\cal X}\subset X$ is a convex subset. 
\end{definition}

In this case, (\ref{pers:utama}) is linear and there are a set of constrainsts (SC) is specified only linear equalities and inequalities, such that there is a polytope if it is bounded. Specifically, the optimization problems form a polyhedron in conditions that (OF)-linear and (SC)-linear. Formally, these as follow \cite{matousek2007,leunberger2008}.

\begin{definition}
[Linear Programming] 
\label{def:lp}
Linear programming (LP) is problems in canonical form:
\begin{equation}
\label{pers:subjeklp}
\max {\bf c}^T{\bf x}
\end{equation}
subject to
\begin{equation}
\label{pers:kendalalp}
A{\bf x}\leq {\bf b}
\end{equation}
where ${\bf x}$ is vector of variables, {\bf c} and {\bf b} are vectors of coefficients, $A$ is a matrix of coefficients.
\end{definition}

\begin{lemma}
\label{lemma:lpco}
If (\ref{pers:subjeklp}) and (\ref{pers:kendalalp}) respect to a linear objective function and constraints are convex, then LP is a convex optima.
\end{lemma}

In case, (\ref{pers:utama}) involves certain types of quadratic programs, we have optimization problems as follow \cite{mittelmann2003,kim2003}.

\begin{definition}
[Second-Order Cone Programming] 
\label{def:socp}
Second-Order Cone Programming (SOCP) is a problem of the form 
\begin{equation}
\label{pers:subjeksocp}
\min f^Tx
\end{equation}
subject to
\begin{equation}
\label{pers:kendalasocp}
\begin{array}{rcl}
\|A_ix+b_i\|_2 &\leq& c_i^T+d_i, i=1,\dots,m\cr 
            Fx &=& g\cr
\end{array}
\end{equation}
where $f\in {\bf R}^n$, $A_i\in {\bf R}^{n_i\times n}$, $b_i\in {\bf R}^{n_i}$, $c_i\in {\bf R}^n$, $d_i\in {\bf R}$, $F\in {\bf R}^{p\times n}$,
$g\in {\bf R}^p$, and $x\in {\bf R}^n$ is the optimization variable.
\end{definition}

\begin{lemma}
\label{lemma:socpco}
If (\ref{pers:subjeksocp}) and (\ref{pers:kendalasocp}) respect to a linear objective function and constraints are convex, then SOCP is a convex optima.
\end{lemma}

\begin{proposition}
\label{proposisi:socplp}
If $A_i=0$, $i=1,\dots,m$ in (\ref{pers:kendalasocp}), then SOCP reduces to LP.
\end{proposition}

Lemma \ref{lemma:lpco}, \ref{lemma:socpco} and Proposition \ref{proposisi:socplp} can be generalized into one conclusion as \emph{semidefinite programming} (SDP)  \cite{kocvara2007,kocvara2009}, where semidefinite matrices as underlying variables, as follow.

\begin{definition}
\label{def:sdp}
Let $S^n$ is the space of all $n\times n$ real symmetric matrices. We define a trace, $tr$, for equipping the space with the inner product, i.e.,
\begin{equation}
\label{pers:sedefinit}
tr(A^TB) = \langle A,B\rangle_{S^n} = \sum_{i=1,j=1}^n A_{ij}B_{ij}.
\end{equation}
A symmmetric matrix is \emph{positive semidefinite} if all its eigenvalues are nonnegative.
\end{definition}

\begin{lemma}
\label{lemma:sdp}
If (\ref{pers:sedefinit}) is convex, then the positive semidefinite is a convex optima.
\end{lemma}

Definition \ref{def:sdp} and Lemma \ref{lemma:sdp} generalize all of forms of convex programming. Thus, LP, SOCP and SDP can be viewed as conic programs with the appropriate type of cone as follow \cite{cezik2005,todd2008,warner2008}.

\begin{definition}
\label{def:co}
Let $X$ is real vector space and also convex. The \emph{conic optimization} is a convex optima with real valued function
\begin{equation}
\label{pers:konikop}
f:C\rightarrow{\bf R}
\end{equation}
defined on convex cone $C\subset X$, and an affine subspace ${\cal H}$ defined by a set of affine constraints $h_i(x) = 0$, such that $x\in C\cap {\cal H}$ for the number $f(x)$ is smallest.
\end{definition}

In geometry, let a monomial is a function $f: {\bf R}^n\rightarrow {\bf R}$ with dom $f={\bf R}_{++}^n$, i.e., $f(x) = cx_a^{a_1}x_2^{a_2}\cdots x_n^{a_n}$, where $c>0$ and $a_i\in {\bf R}$. Based on this concept, we define the geometric programming as follow \cite{rajgopal2002,boyd2007,yang2010}.

\begin{definition} [Geometric programming]
\label{def:gp}
 A geometric programming (GP) is optimization problem of the form
\begin{equation}
\label{pers:subjekgp}
f_0(x)
\end{equation}
subject to
\begin{equation}
\label{pers:kendalagp}
\begin{array}{rcl}
f_i(x)&\leq& 1, i = 1,\dots,m\cr
h_i(x)&=& 1, i=1,\dots,p\cr
\end{array}
\end{equation}
where $f_0,\dots,f_m$ are polynomial and $h_1,\dots,h_p$ monomials.
\end{definition}

\begin{lemma}
\label{lemma:gp}
GP is a convex optima, if (\ref{pers:subjekgp}) and (\ref{pers:kendalagp}) 
become a sum of exponentials of affine functions.
\end{lemma}

Some concepts of optimization problems above include into convex programming, and then we will define the concepts of other programming. One of them is quadratic programming as a special type of mathematical optimization problem \cite{kim2003}.

\begin{definition} [Quadratic Programming]
\label{def:qp}
Let ${\bf x} \in {\bf R}^n$. The $n\times n$ matrix $Q$ is symmetric, and ${\bf c}$ is any $n\times 1$ vector. Quadratic programming (QP) is to minimize (with reppect to {\bf x})
\begin{equation}
\label{pers:subjekqp}
f({\bf x}) = \frac{1}{2}{\bf x}^TQ{\bf x}+{\bf c}^T{\bf x}
\end{equation}
subject to one or more constraints of of the form
\begin{equation}
\label{pers:kendalaqp}
\begin{array}{rcl}
A{\bf x}&\leq& {\bf b}\cr
E{\bf x}&=&{\bf d}\cr
\end{array}
\end{equation}
where ${\bf x}^T$ indicates the vector transpose of ${\bf x}$. 
\end{definition}

The notation $Ax\leq b$ means that every entry of the vector $Ax$ is less than or equal to the corresponding entry of the vector ${\bf b}$.

\begin{lemma}
If the matrix $Q$ is positive semidefinite, then (\ref{pers:subjekqp}) is a convex function, or QP is convex optima.
\end{lemma}

\begin{definition} [Nonlinear Programming] 
\label{def:nlp}
\emph{\cite{leunberger2008}}
Optimization model as nonlinear programing is problem can be stated simply as
\begin{equation}
\label{pers:nonlinear1}
\max_{x\in X} f(x)
\end{equation}
is to maximize some variable such as product throughput, or
\begin{equation}
\label{pers:nonlinear2}
\min_{x\in X} f(x)
\end{equation}
is to minimize a cost function, where $f:R^n\rightarrow R$, $X\subset R^n$.
\end{definition}

All models above are the determenistic optimization problems which are formulated with known parameters. Besides, the size of $X$ is bounded, also variable with another have self setting. When the parameters are known within certain bounds, one approach to tackling such problem is called \emph{robust optimization}. However the real world problems almost invariably include some unknown parameters. The goal of optimization is to find a solution which is feasible for all conditions and situations and optimal in some sense. Therefore, one of frameworks for modeling optimization problems is by involving uncertainty. 

\begin{table}
\caption{Methods are used to solve optimization problems.}
\label{tabel:metode}
\begin{center}
\begin{tabular}{|r|l|l|}\hline
id & Method & Description\cr\hline
1. & active set \cite{murty1998} & A problem is defined using an objective function to\cr 
   &                             & minimize or maximize, and a set of constraints $g_1(x)$\cr
   &                    & $\geq 0,\dots,g_k(x)\geq 0$ that define the feasible region, that \cr   
   &                             & is, the set of all x to search for the optimal solution.\cr
   &                             & Given a point x in the feasible region, a constraint\cr
   &                             & $g_i(x)\geq 0$ is called active at $x$ if $g_i(x) = 0$ and \cr
   &                             & inactive at $x$ if $g_i(x) > 0$. \cr
2. & ant colony \cite{colorni1991} & A probabilistic technique for solving computational\cr
   &                             & problems which can be reduced to finding good paths \cr
   &                             & through graphs.\cr
3. & beam search \cite{kim2004,bautista2008} & A heuristic search algorithm that explores a graph by \cr
   &                             & expanding the most promising node in a limited set.\cr
4.   & conjugate gradient          & An algorithm for the numerical solution of particular \cr
   & \cite{dai2001}                & systems of linear equations, namely those whose  \cr
   &                             & matrix is symmetric and positive-definite\cr
5.   & cuckoo search \cite{yang2009} & An optimization algorithm was inspired by the  \cr
   &                             & obligate brood parasitism of some cuckoo species by \cr
   &                             & laying their  eggs in the nests of other host birds \cr
   &                             & (of other species). \cr
6.   & differential evolution      & A method that optimizes a problem by iteratively \cr
   &  \cite{zhang2008,zhang2009}  & trying to improve a candidate solution with regard to \cr 
   &                             & a given measure of quality.\cr
7.   & dynamic relaxation & A numerical method, which, among other things, can  \cr
   &   \cite{lawphongpanich2006} & be used do "form-finding" for cable and fabric \cr
   &                             & structures.\cr
8.   & ellipsoid \cite{grotschel1981} & An iterative method for minimizing convex functions.\cr
9.   & evolution strategy          & An optimization technique based on ideas of  \cr
   & \cite{gonzalez1998}          & adaptation and evolution. \crcr
10.   & firefly algorithm           & A metaheuristic algorithm, inspired by the flashing \cr
   &  \cite{lukasik2009,yang2010b}    & behaviour of fireflies.\cr
 11.  & Frank-Wolfe \cite{chryssoverghi1997}& A procedure for solving quadratic programming  \cr
   &                             & problems with linear constraints.\cr
 12.  & genetic algorithm           & A search heuristic that mimics the process of natural\cr
   & \cite{onwubolu2003,dogan2004}            & evolution. This heuristic is routinely used to generate\cr
   &                             & useful solutions to optimization and search problems.\cr
 13.  & gradient projection         & An algorithms that can be used to optimize virtually \cr
   &  \cite{guocheng2010}        & any rotation criterion.\cr
 14.  & harmony search              & A phenomenon-mimicking algorithm inspired by the\cr
   &  \cite{rappaport2007,kaveh2009}           & improvisation process of musicians for finding \cr
   &                             & a best harmony as global optimum.\cr
 15.  & hill climbing               & A mathematical optimization technique which belongs \cr
   &  \cite{greiner1996,lewis2009}            & to the family of local search.\cr
 16.  & interior point  &  A linear or nonlinear programming method that\cr
   &  \cite{forsgren2002,bonnans2006}  & achieves optimization by going through the middle of \cr
   &                             & the solid defined by the problem rather than around \cr
   &                             & its surface. \cr\hline
\end{tabular}
\end{center}
\end{table}

Ontologically, in certain optimization problems the unknown optimal solution might not be a number or vector, but rather a continuous quantity \cite{jeyakumar2008,pham2009}. This problem appeared because a continuous quantity cannot be determined by a finite number of certain degrees of freedom \cite{bot2009,chen2005,mishra2008}. However, such problem can be more challenging than finite-dimensional ones \cite{hinze2009}, so nothing of entities attributes depend on continuous times ever, or there are not event hold on long times. We know disciplines which study infinite-dimensional optimization problems are \emph{calculus of variations} \cite{buttazzo2009,chen2005}, \emph{optimal control} \cite{chinchuluun2010} and \emph{shape optimization} \cite{bucur2005,kovtumenko2006,eppler2008}. 
In constraints based optimization \cite{hinze2009}, a solution is a vector of variables that satisfies all constraints, where the constraint satisfication is the process of finding a solution to a set of constraints that impose conditions that the variables must satisfy. However, the solutions depend on Adaption 1, where constraint satisfaction typically identified with problems based on constraints on a finite domain, or based on patterns derived from the experience in classification, that be local solution.
\begin{table}
\caption{Methods are used to solve optimization problems (Continue)}
\label{tabel:metode2}
\begin{center}
\begin{tabular}{|r|l|l|}\hline
id & Method & Description\cr\hline
 17.  & IOSO                        & A multiobjective, multidimensional nonlinear \cr
   &                             & optimization  technology.\cr
 18.  & line search  \cite{belegundu2004,shi2005}                & One of two basic iterative approaches to finding a\cr
   &   & local minimum ${\bf x}^*$ of an objective function $f:{\bf R}^n\rightarrow {\bf R}$.\cr
 19.  & Nelder-Mead  \cite{price2002,burmen2005}               & A nonlinear optimization technique, which is a \cr
   &     & well-defined numerical method for twice differentiable \cr
   &                             & and unimodal problems.\cr
 20.  & Newton \cite{mangasarian2004,fujishige2009}                     & A method for finding successively better approximations\cr
   &                             & to the zeroes (or roots) of a real-valued function.\cr
 21. & particle swarm \cite{wu2009} & A computational method that optimizes a problem by\cr     
   &                             & iteratively trying to improve a candidate solution \cr
   &  	                         & with regard to a given measure of quality.\cr
 22.  & quantum annealing           & A general method for finding the global minimum of a\cr
   &                             & given objective function over a given set of candidate\cr     
   &                             & solutions (the search space), by a process analogous \cr
   &                             & to quantum fluctuations.\cr
23.   & quasi-Newton \cite{jiang2003} & A algorithm for finding local maxima and minima of \cr
   &                             & functions.\cr
24. & simplex \cite{ehrgott2007}  & A popular algorithm for numerically solving linear\cr
   &                             & programming problems. \cr
25.  & simulated annealing         & A generic probabilistic metaheuristic for the global\cr
   &  \cite{yen2004}             & optimization problem of applied mathematics, namely \cr
   &                             & locating a good approximation to the global optimum of a\cr
   &                             & given function in a large search space.\cr
26.  & stochastic tunneling     & An approach to global optimization based on the Monte \cr
   & \cite{oblow2001}           & Carlo method-sampling of the function to be minimized.\cr
27.   & subgradient                 & The iterative methods for solving convex minimization \cr
   & \cite{gokbayrak2009}        & problems.\cr
28.  & tabu search                 & A mathematical optimization method for local search \cr
   & \cite{nowicki2005}          & techniques.\cr\hline
\end{tabular}
\end{center}
\end{table}

In stochastic framework, some of the constraints or parameters depend on random variables \cite{pham2009,schneider2006,zhigljavsky2008}, or this models of optimization depend on probability distributions governing the data are known or can be estimated. The goal is to find some policy that is feasible for all (or almost all) the possible data instance and maximizes the expectation of some function of the decisions and the random variables. Although, stochastic programming has applications in a broad range of areas (most of them in uncertainty): finance, transportation, social modalities, or energy. However, random event occurs only affecting first acitivity.

In theoretical computer science, many themes of optimization involving combinatorial \cite{avis2005,du2005,iglesias2005}, that related to operations research, algorithm theory, and computational complexity theory. 
The goal is to find the best solution. However, this optimization problems only dealing with graphs, matroids, and related structures, in which the set of feasible solutions is discrete or can be reduced to discrete \cite{kasperski2008,lozovanu2009,papadopulos2009}. In addition to, for combinatorial optimization are designed metaheuristic that optimizes a problem by iteratively trying to improve a candidate solution with regard to a given measure of quality \cite{doerner2007,rego2005,siarry2008}, where an optimal solution is sough over a discrete search-space. However, metaheurstic do not guarantee an optimal solution is ever found. Therefore, many metaheuristics implementation in stochastic optimization \cite{sorensen2003,tseng2006,lewis2007}.

Many models of optimization have been defined, and many methods for solving them have been created to what(?). In modelling optimization, the existence of derivatives is not always assumed, therefore many methods were devised for specific situations. All methods are created based on smoothness of the objective functions (\ref{pers:utama}), like as class of combinatorial methods, class of derivative-free methods, class of first-order methods, class of second-order methods, see Tabel \ref{tabel:metode} and \ref{tabel:metode2}. However, all methods lie on convexity \cite{burachik2008}: should the objective function be convex over the region of interest, then any local minimum will also be a global minimum. So, the ontology of optimization models and their methods is used to describe chain one optimization problem with another into a structure as knowledge so that knowledge based optimization exists, i.e. a systematic approach for the study of optimization indiscernibility, where the possibilities of attributes (as independent variables) will be playing \cite{nasution2010}.

\section{Discussion}

All models of optimization (Definition \ref{def:lp}-\ref{def:nlp}) came closer to Adaptation 2, and methods for solving them to. Few methods have shown the possibility of solution using Adaptation 3, among of them in implicit conditions. Of course rasionally, many models of optimization have been defined in one space ${\bf R}^n$, where a function employed to relate between one (may be two) variable to many variables. Ontologically, the real world consist of relations between multiple spaces with variety of dimensions. For example, each academic persons related to attributes: name, age, salary, papers, etc. Each paper has a set of descriptions as attribute, i.e., title, venue, co-author, publisher, pages, year, etc. So, a publisher as entity related to: name, city, country, etc. There is an attribute has discrete type (name of person), but has relation with another atribut in different space (name of co-author) in probability, or there are a stochastic relation between name of person and name of co-author laying (social network) on title of paper by using Bayesian approach \cite{nasution2010b}. 

The optimization problems not only about looking for answers of questions:
\begin{enumerate}
\item [$~~~~~~~~~~~~~~~~~~~~~~~~~~~~~~$] Is it possible to satisfy all constrains?
\item [$~~~~~~~~~~~~~~~~~~~~~~~~~~~~~~$] Does an optimum exists?
\item [$~~~~~~~~~~~~~~~~~~~~~~~~~~~~~~$] How can an optimum be found?
\item [$~~~~~~~~~~~~~~~~~~~~~~~~~~~~~~$] How does the optimum change if the problem changes?
\end{enumerate}
In the opposite, about question:

\begin{center}
(MQ)\\
 {\it Why are so many models and methods of optimization created and so few used}? 
\end{center}

We start this discussion using one analogy: The optimization problems are often expressed by $\min_{x\in {\bf R}} f(x)$ or $\max_{x\in{\bf R}} f(x)$. Certainly, $x$ ranges over the reals. What is answer for ${\rm argmin}_{x\in(-\infty,-1]} f(x)$. This asks for one (or more) values of $x$ in the interval $(-\infty,-1]$ that minimize the function. The answer is the undefined. This is satisfiability problem, i.e., the problem of finding any feasible solution at all without regard to objective value (do not satisfy (MQ)). This means that any feasible solution do not realiable! Why? Many optimization method need to start from a feasible point. What? One way to obtain such a point is to relax the feasibility conditions using a slack variable, with enough slack, any starting point is feasible. How? To minimize that slack variable until slack is null or negative. So, where is the positition of optimization goal?

\begin{theorem}
[Extreme Value of Karl Weierstrass] If a real-valued function $f$ is continuous in the closed and bounded interval $[a,b]$, then $f$ must attain its maximum and minimum value, each at least once.
\end{theorem}

\begin{table}
\caption{Description of linear programming}
\label{des:lp1}
\begin{center}
\begin{tabular}{|rcl|}\hline
Name && {\bf Linier Programming} [Convex Optima]\cr\hline
Description && A mathematical method for determining a way to achieve the best \cr
            && outcome in a given mathematical model for some list of requirements\cr
     && represented as linear equations.\cr\hline
Formula && {\it Canonical form}: Definition \ref{def:lp}\cr
        && Standard form:\cr 
        && 1. A linear function to be maximized\cr
        && ~~~e.g., maximize ${\bf c}_1{\bf x}_1+{\bf c}_2{\bf x}_2$\cr
        && 2. Problem constraints of the following form, i.e., \cr
        && ~~~$a_{1,1}{\bf x}_1+a_{1,2}{\bf x}_2\leq {\bf b}_1$\cr
        && ~~~$a_{2,1}{\bf x}_1+a_{2,2}{\bf x}_2\leq {\bf b}_2$\cr
        && ~~~$a_{3,1}{\bf x}_1+a_{3,2}{\bf x}_2\leq {\bf b}_3$\cr
        && 3. Non-negative variables, e.g., ${\bf x}_1\geq 0$, ${\bf x_2}\geq 0$.\cr
        && 4. Non-negative right hand side constants ${\bf b}_i\geq 0$.\cr\hline
\end{tabular}
\end{center}
\end{table}
This theorem describes optimum solution exist under conditions, explicitly, i.e., there exist numbers $c$ and $d$ in $[a,b]$ such that: $f(c)\leq f(x)\leq f(d)$, $\forall x\in [a,b]$.

\begin{theorem}
[Fermat theorem] Let $f:(a,b)\rightarrow{\bf R}$ be a function and suppose that $x_0\in(a,b)$ is a local extremum of $f$. If $f$ is differentiable at $x_0$ then $f'(x_0) = 0$.
\end{theorem}

Last theorem states that optima of unconstrained problems are found at stationary points, where the first derivative or the gradient of the objective function is zero. In general, the optima may be found at critical points, where the first derivate or gradient of the function is zero or is undefined, or on the boundary of the choosed set. What? First derivative test only identifies points that might be optima, but it cannot distinguish a point which is a minimum from one that is a maximum or one that is neither. How? Do twice differentiable, and check the second derivatives (matrix) for unconstrained problems. Therefore, the conditions that distinguish maxima and minima from other stationar points called the second-order conditions. 
While the inequality-constrained optimazation problems are conditioned by the Lagrange multiplier, or calculating the complementary slackness conditions based on Karush-Kuhn-Tucker conditions. Then, check the matrix of second derivatives of the objective function and the constraints. So, how many steps to find optimum solution in order to it is applicable?  

The convexity of (\ref{pers:convexoptima}) in Definition \ref{def:ko} makes the powerful tools of convex analysis applicable, such as theory of subgradients lead to a particularly satisfying theory of necessary and sufficient conditions for optimality, based on theorem as follows:

\begin{theorem}
[Hahn-Banach] Given a vector space $V$ over the field ${\bf R}$ of real numbers, a function $f:V\rightarrow {\bf R}$ is called sublinear if $f(\gamma x) = \gamma f(x)$ for any $\gamma\in {\bf R}_+$ and any $x\in V$, $f(x+y)\leq f(x)+f(y)$ for any $x,y\in V$.
\end{theorem}

In ontology, the applicability of linear programming is proved starting in Table \ref{des:lp1}. In general, a standard form is usual and most intuitive form of describing a convex minimization problem: (1) A convex function $f(x) : {\bf R}^n\rightarrow {\bf R}$ to be minimized over the variable $x$; (2) The constraints, if any: (a) Inequality constraints of the form $g_i(x)\leq 0$, where the functions $g_i$ are convex. (b) Equality constraints of the form $h_i(x)=0$, where the function $h_i$ are linear. 

\begin{table}
\caption{Description of linear programming}
\label{des:lp2}
\begin{center}
\begin{tabular}{|rcl|}\hline
Name && {\bf Linier Programming} [Convex Optima]\cr\hline
Description && A type of convex programming, a model optimiztion where objective \cr
            && is affine and the set of constraints is in affine.\cr\hline
Formula && 1. Primal problem:\cr
        && ~~~Maximize ${\bf c}^T{\bf x}$ subject to $A{\bf x}\le {\bf b}$, ${\bf x}\geq 0$.\cr
        && 2. Symmetric dual problem:\cr
        && ~~~Minimize ${\bf b}^T{\bf y}$ subject to $A^T{\bf y}\geq {\bf c}$, ${\bf y}\geq 0$.\cr
        && 3. Alternative primal formulation:\cr
        && ~~~Maximize ${\bf c}^T{\bf x}$ subject to $A{\bf x}\le {\bf b}$.\cr
        && 4. Symmetric dual problem:\cr
        && ~~~Minimize ${\bf b}^T{\bf y}$ subject to $A^T{\bf y} = {\bf c}$, ${\bf y}\geq 0$, ${\bf y}\geq 0$. \cr\hline
\end{tabular}
\end{center}
\end{table} 

In optimization model, the duality theory generalize convex programming where methods took an effective computational. Convex programming minimizes convex functions, so it useful also for maximizing concave functions, see Table \ref{des:lp1} and Table \ref{des:lp2}: The problem of maximizing a concave function can be reformulated equivalently as a problem of minimizing a convex function. 
Thus, a convex minimization problem written as $\min_x f(x)$ subject to $g_i(x)\leq 0$, $i=1,\dots,m$, $h_i(x) = 0$, $i=1,\dots,p$, where every equality constraint $h(x) = 0$ can be equivalently replaced by a pair of inequality constraint $h(x) \leq 0$ and $-h(x)\leq 0$. By this, $h_i(x) = 0$ has to be affine as opposed to merely being convex. If $h_i(x)$ is convex, $h_i(x)\leq 0$ is convex, but $-h_i(x)\leq 0$ is concave. Therefore, the only way for $h_i(x) = 0$ to be convex is for $h_i(x)$ to be affine. The way is a transformation!

Can we see that the optimization problems is just one problem only? What is answer! For linear programming, we have a solution with many methods, see Table \ref{des:lp3}. 

\begin{table}
\caption{Description of linear programming}
\label{des:lp3}
\begin{center}
\begin{tabular}{|rcl|}\hline
Name && {\bf Linier Programming} [Convex Optima]\cr\hline
Description && A type of convex programming, a model optimiztion where objective \cr
            && is affine and the set of constraints is in affine.\cr\hline
Formula && Canonic form $\models$\cr
        && ~~~~~~~~~~1. Primal problem\cr
        && ~~~~~~~~~~2. Dual problem\cr\hline
\end{tabular}
\begin{tabular}{|l|l|l|l|}\hline
        & Founder            & Year & Solution Technique\cr \hline
Model & Leonid Kantorovich & 1939 & \cr
Methods & George B. Dantzig  & 1947 & Simplex Method\cr
        & John von Neumann   & 1947 & Duality\cr
        & Leonid Khachiyan   & 1979 & Polynomial Time\cr
        & Narendra Karmarkar & 1984 & Interior point method\cr\hline
\end{tabular}
\end{center}
\end{table} 

Consider an arbitrary maximization (or minimization) problem where the objective function $f({\bf x},{\bf r})$ depends on some parameters ${\bf r}$, $f^*({\bf r}) = \max_{\bf x} f({\bf x},{\bf r})$. The function $f^*({\bf r})$ is the problem's optimal-value function gives the maximized or minimized value of the objective function $f({\bf x},{\bf r})$ as a function of its parameter ${\bf r}$.
 
\begin{theorem}
[Envelope theorem] Let $\max_{\bf x} f({\bf x},{\bf r})$ s.t. $g({\bf x},{\bf r}) = 0$, and Lagrangian function ${\cal L}({\bf x},{\bf r}) =  f({\bf x},{\bf r}) - \lambda g({\bf x},{\bf r})$, where $\lambda = (\lambda_1,\dots,\lambda_n)$, $g({\bf x},{\bf r}) = (g_1({\bf x},{\bf r}),\dots,g_n({\bf x},{\bf r}))$, ${\bf 0} = (0,\dots,0)\in {\bf R}^n$. Then 
\begin{equation}
\frac{\partial f^*({\bf r})}{\partial r_i} = \frac{\partial{\cal L}({\bf x},{\bf r})}{\partial r_i}\Big|_{{\bf x}={\bf x}*({\bf r}),\lambda = \lambda({\bf r})}.
\end{equation}
\end{theorem}

This theorem expresses how the value of an optimal solution changes when an underlying parameter changes. While, the continuity of the optimal solution as a function of underlying parameters, i.e. 

\begin{theorem}
[Maximum theorem] Let $X$ and $\Theta$ be metric space, $f:X\times\Theta\rightarrow {\bf R}$ be a function jointly continuous in its two arguments, and ${\cal C}:\Theta\rightarrow X$ be a compact-valued correspondance. For $x\in X$ and $\theta\in \Theta$, let $f^*(\theta) = \max\{f(x,\theta)|x\in C(\theta)\}$ and $C^*(\theta) = {\rm arg}\max\{f(x,\theta)|x\in C(\theta)\} = \{x\in C(\theta)|f(x,$ $\theta) = f^*(\theta)\}$. If $C$ is continuous at some $\theta$, then $f^*$ is continuous at $\theta$ and $C^*$ is non-empty, compact-valued, and upper at $\theta$.
\end{theorem}

Consider the restriction of a convex function to a compact convex set, where the function on that set attains its constrained maximum only on the boundary. For example, Definition \ref{def:gp} defines that objective function and inequatlity constrainsts of GP can be expressed as \emph{psynomials} and equality constrainsts as \emph{monomials}. Both psynomials and monomials can be transformed into a convex program (Lemma \ref{lemma:gp}). Other side, constrained problems can often be transformed into unconstrained problems by using Lagrange multipliers.

Consider a convex minimization problem given in standard form by a cost function $f(x)$ and inequality constraints $g_i(x)\leq 0$, where $i=1,\dots,m$. Then the domain ${\cal X}$ is ${\cal X} = \{x\in X|g_1(x)\leq 0,\dots,g_m(x)\leq 0\}$. The Lagrangian function for this problem is $L(x,\lambda_0,\dots,\lambda_m) =\lambda_0f(x)+\lambda_1g_1(x)+\dots+\lambda_mg_m(x)$. For each point $x \in X$ that minimizes (\ref{pers:utama}), there exist real number $\lambda_0,\dots,\lambda_m$, called Lagrangian multipliers, that satisfy these conditions simultaneously, (a) $x$ minimizes $L(y,\lambda_0,\dots,\lambda_m)$ over all $y\in X$, (b) $\lambda_0\geq 0$, $\lambda_1\geq 0$, $\dots$, $\lambda_m\geq 0$, with at least one $\lambda_k>0$, and (c) $\lambda_1g_1(x) = 0$, $\dots$, $\lambda_mg_m(x) = 0$. For instance, if there is a strictly feasible point, namely a point $z$ satisfying $g_1(z)<0,\dots, g_m(z)<0$, then must be assigned $\lambda_0=1$. Conversely, if some $x\in X$ satisfies (a)-(c) for $\lambda_0,\dots,\lambda_m$ with $\lambda_0=1$, then $x$ is certain to minimize (\ref{pers:utama}). Therefore,
\begin{lemma}
Three following statements is equivalent.
\begin{enumerate}
\item If a local minimum exists, then its is a global minimum.
\item The set of all minima is convex.
\item For each strictly convex function, if the function has a minimum, then the minimum is unique.
\end{enumerate}
\end{lemma}

\begin{proposition}
A set of problems that consist of least squares, LP, QP, Conic Optimization, GP, SOCP, SDP, Quadratically Constrained Quadratic Programming, and Entropy Maximization are convex optima if and only if can be transformed into convex minimization problems by changing variables.
\end{proposition}

The last proposition is typically interpreted as prividing conditions for optimization problems, where their ontology describe as follow:
\[
\begin{array}{rclcccc}
{\rm Optimization~problems} &\models& {\rm Convex~optima}     &\dashv  & &\cr
                            &\models& {\rm Conic~Programming} &\uparrow&\dashv& \cr
                            &\models& {\rm LP}                &{\rm Lemma~1}&\uparrow&\dashv\cr
                            &\models& {\rm SOCP}              &{\rm Lemma~2}&\uparrow&{\rm Prop.~1}\cr
                            &\models& {\rm SDP}               &{\rm Lemma~3}&\perp&\cr
                            &\models& {\rm GP}                &{\rm Lemma~4}& &\cr
                            &\models& {\rm QP}                &{\rm Lemma~5}& &\cr
                            &\models& {\rm Non-LP}                & & &\cr
                            &\models& {\rm Stochastic~Programming}& & &\cr
                            &\models& {\rm Fuzzy~Programming}& & &\cr
                            &\models& {\rm ?~Programming}& & &\cr
{\rm Opt.~Indiscernebility}&\models& {\rm Incomplete~Programming}& & &\cr
\end{array}
\]

To complete the interpretation of optimization problems and to understand the incomplete with them, some considerations have formulated as advice \cite{landry1996}, i.e., (a) the model of optimization be ready in close with strategic stakeholder with understanding the organization or environment; (b) the model of optimization is designed to adapt task at hand and to the cognitive capacity of the stakeholder; (c) the model become familiar with the various logics and preferences prevailing in the organization or many other environments; (d) the model of optimization problems must be formulated in comportable and a framework where many methods can be solved them; (e) the model optimization problems can be prepared to modify or develop a new version, change the framework in order to accordance with any method; (f) each model can give the options to select the making decisions; (g) the models can express the respecting implicitly become explicit.

\begin{conjecture}
The optimization in incomplete has optimum solution.
\end{conjecture}

\section{Conclusion}
Throughout the ontology of optimization problems were in close connection with the principles and findings of the science, philosophy, and arts. It can be noticed that optimization problems are a reflection of pictures respecting to entities, where the ontology is the beginning of philosophy, the scratch of art, the formalization in science, to generalize knowledge behind a model of optimization.

\begin{remark}
The complete is dual of incomplete.
\end{remark}

\end{document}